\documentclass[english]{article}
\usepackage[T1]{fontenc}
\usepackage[latin9]{inputenc}
\usepackage[letterpaper]{geometry}
\usepackage{amsmath}
\usepackage{authblk}
\usepackage{graphicx}
\usepackage{epstopdf}
\usepackage{color}
\usepackage{babel}

\usepackage{amsfonts} 

\usepackage[notref,notcite]{showkeys}

\graphicspath{ {./results-perforated/}}

\makeatletter

\title{Re-iterated multiscale model reduction using the GMsFEM}

\author[1,2]{
Eric T. Chung \thanks{Email: {\tt tschung@math.cuhk.edu.hk}.} }
\author[2]{
Yalchin Efendiev \thanks{
 Email: {\tt efendiev@math.tamu.edu}.}}
\author[2]{
Wing Tat Leung}
\author[2,3]{
 Maria Vasilyeva}

\affil[1]{Department of Mathematics,
The Chinese University of Hong Kong (CUHK), Hong Kong SAR. }

\affil[2]{Department of Mathematics \& Institute for Scientific Computation (ISC),
Texas A\&M University,
College Station, TX 77843-3368, USA.}

\affil[3]{Department of Computational Technologies, Institute of Mathematics and Informatics, North-Eastern Federal University, Yakutsk, 677980, Republic of Sakha (Yakutia), Russia. }

\begin{document}
\maketitle
%===========

\begin{abstract}

Numerical homogenization and multiscale finite element methods
construct effective properties on a coarse grid by solving
local problems and extracting the average effective properties
from these local solutions.
In some cases, the solutions of local problems can be expensive to
compute
due to scale disparity. In this setting, one can basically apply
a homogenization or multiscale method re-iteratively
to solve for the local problems. This process is known as re-iterated
homogenization and has many variations in the numerical context.
Though the process seems to be a straightforward extension of
two-level process, it requires some careful implementation and the
concept development for problems without scale separation and
high contrast. In this paper, we consider the Generalized
Multiscale Finite Element Method (GMsFEM) and apply it iteratively to
construct its multiscale basis functions. The main idea of the GMsFEM
is to construct snapshot functions and then extract multiscale basis
functions (called offline space)
 using local spectral decompositions in the snapshot spaces.
The extension of this construction to several levels  uses snapshots
and offline spaces interchangebly to achieve this goal.
At each coarse-grid scale,
we assume that the offline space is a good approximation
of the solution and use all possible offline functions or
randomization as boundary conditions and solve the local problems
in the offline space at the previous (finer) level, to construct
snapshot space. We present an adaptivity strategy and show numerical
results for flows in heterogeneous media and in perforated domains.

\end{abstract}

\section{Introduction}

%General about multiscale and the need for multilevel
Many problems with multiple scales can have a large scale disparity.
For example, porous media problems can have spatial variations from the pore-scale
sizes to the field
scales. Solving these problems requires coarsening
approaches. Some successful coarsening methods include homogenization
\cite{blp78,ab05, jko94},
numerical homogenization and generalization \cite{weh02,dur91,fish_book,eh09,fan2010adaptive, fish2012homogenization,li2008generalized}, and multiscale methods
\cite{fish_book,ee03,abe07,oh05,hw97,eh09, ehw99, egw10,fish2012staggered,chung2016_review,efendiev2014generalized,efendievsparse, franca2005multiscale,nouy2004multiscale,chung2016mixed,calo2016multiscale,gao2016application,chung2013sub}.
The objective of these approaches is to solve the problem on a prescribed
computational grid, which we will call the coarse grid. The coarse grid can typically be
much larger than the fine grid, which resolves the underlying processes.
These coarse-grid techniques require local computations to extract
effective properties. In numerical homogenization methods, this involves
solving cell problems. In multiscale methods, this involves
solving multiscale basis functions, where multiscale
basis functions play a role of effective properties
\cite{chung2016_review}.
 In problems with very large disparate scales,
the solution of local problem can be expensive and require
coarsening.

%Existing multilevel methods and they are advantages. A. Brandt. What is different
%in general words.
There have been several approaches to approximate the local solutions, when
there is very large scale disparity.
The commonly used technique, which motivates our studies, is
re-iterated homogenization, which goes
back to  \cite{blp78}, see also
\cite{cioranescu2002periodic,lions2001reiterated}.
Other methods include
hierarchical basis functions  \cite{hoang2005high,schwab2003high},
hierarchical model reduction \cite{yuan2009hierarchical,fish_book,fish1997hierarchical},
multilevel approaches for multiscale basis
functions \cite{lipnikov2008multilevel,maclachlan2006multilevel,lipnikov2011adaptive},
hierarchical upscaling \cite{brown2013efficient,fish1993multiscale,oden1999hierarchical},
multilevel approximations \cite{brandt91,brandt99,MLBFP,vassilevski_upscaling,egv2011},
 and approximate basis computations
 \cite{eberhard2004coarse}.
In these approaches, multiscale basis functions
or local solutions are approximated in a re-iterated fashion by constructing
multiple level approximations for the basis functions.
The accuracy of local solutions may or may not have a large effect on
the accuracy of the coarse-grid solutions due to subgrid errors dependent on
how accurate localized approaches.
These approximate local techniques can demonstrate an efficiency and
speed-up in practical computations in many situations.
However, these multiscale approaches construct a limited number of basis functions in
each coarse block and mainly focus on designing an approximation
step for calculating pre-defined local problems.
Our objective is to extend a systematic approach
proposed in \cite{egh12,chung2016_review} in a re-iterated fashion and
show that one can provide a systematic
way of eliminating the degrees of freedom using snapshots and local spectral problems.
%Our approach shares common concepts with the multilevel approximation techniques
%a common concept with multilevel upscaling proposed in
% \cite{brandt91,brandt99,MLBFP,vassilevski_upscaling,egv2011}.

%The method. The main idea
To describe the main idea of the proposed method, we briefly mention the underlying
concept of the GMsFEM.
The main idea of the GMsFEM is to introduce snapshots on a coarse grid and
identify dominant modes in the snapshot space. The snapshot space represents
a set of functions that can be used to
accurately
 calculate the local solution space.
The snapshots typically consist of local solutions with all possible
 boundary conditions
or with randomized boundary conditions, which
represent point sources distributed on the boundaries.
The snapshot vectors are similar to the snapshots used in
global model reduction \cite{chung2016_review}; however,
they are constructed without solving expensive
global problems and similar to cell problems in homogenization.
 The snapshot vectors are computed either by using all boundary conditions
or by using random boundary conditions
\cite{chung2016_review}.
To avoid effects due to oscillations of random boundary conditions
and improve the accuracy, we use oversampling techniques and compute
snapshots in domains slightly larger than the target coarse block.
The local spectral decomposition is performed
in the snapshot space by identifying local spectral decomposition based
on the analysis. The analysis consists of decomposing the error into local
regions and bounding these error components. Furthermore, adaptivity is used
to identify the regions that require more basis functions.

\begin{figure}[!htb]
  \centering
  \includegraphics[width=1\textwidth]{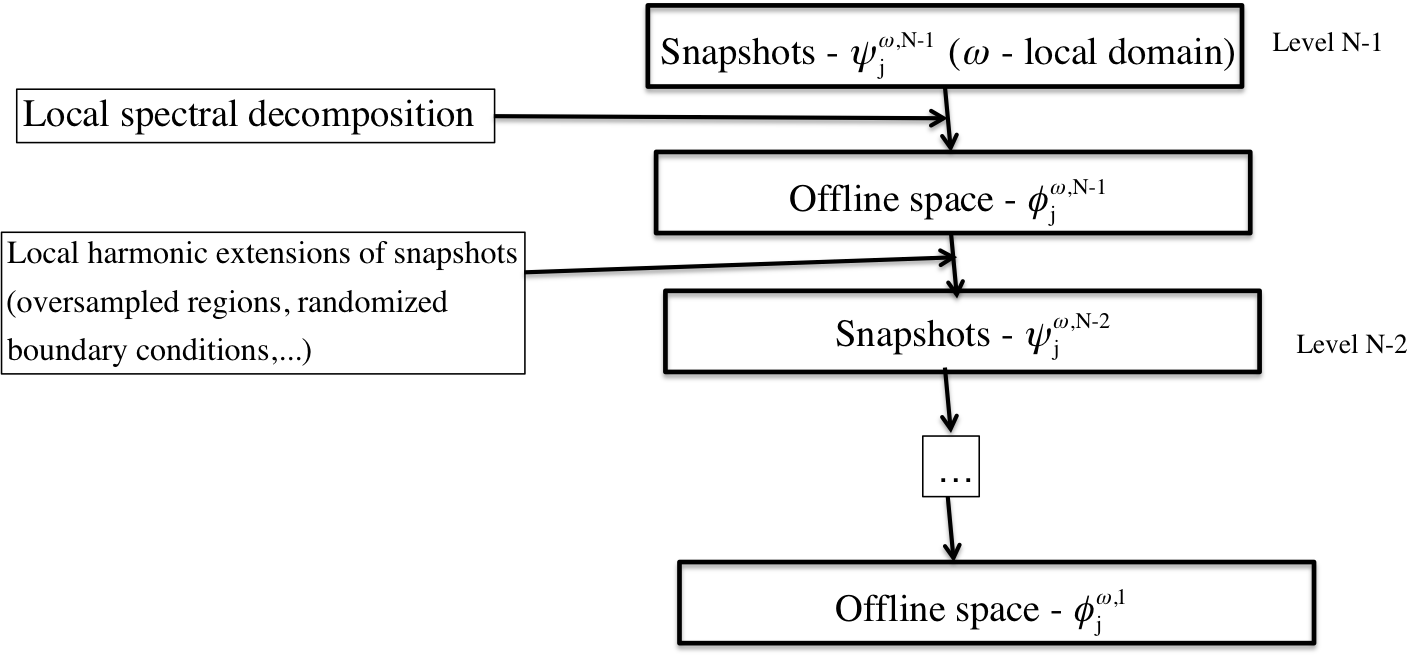}
  \caption{Illustration. The outline of the algorithm.}
 \label{Illustration}
\end{figure}

\begin{figure}[!htb]
  \centering
  \includegraphics[width=1\textwidth]{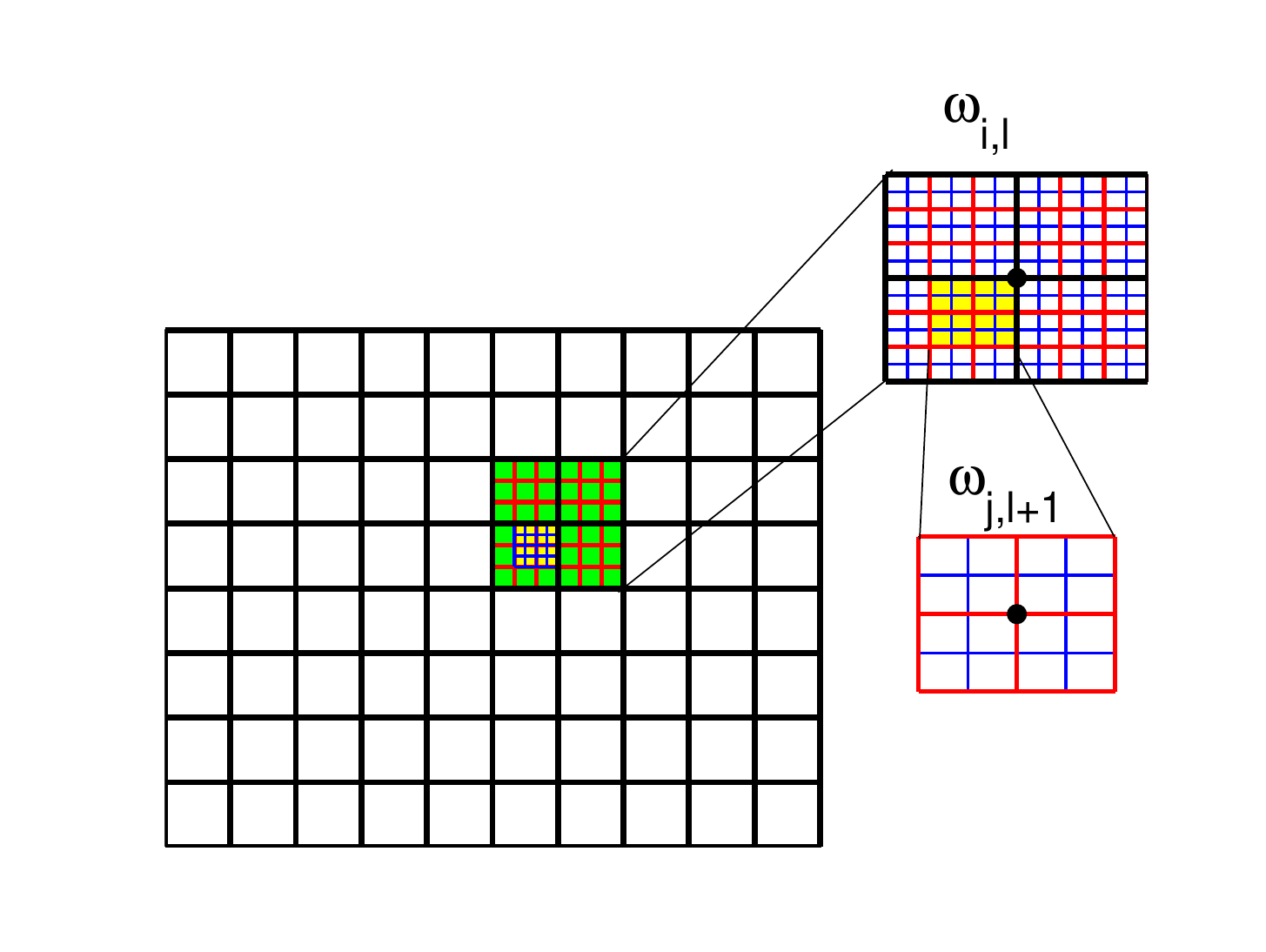}
  \caption{Mesh}
 \label{Mesh}
\end{figure}

%The ingredients
Extending
this approach to several disparate scale setups requires repeating the above
process. Assume that we have a hierarchy of coarse meshes $H_1>H_2>...>H_N=h$,
where our objective is to compute multiscale basis functions on the coarse grid $H_1$
(see Figure \ref{Mesh} for illustration).
Here, for simplicity, $H_l$ denotes the mesh size,
also referred to
% as well as we refer it
as a mesh
configuration. The main contribution in extending two-level approaches
 to the multi-level
is to identify snapshots at each level.
We depict an illustration of the main concepts in Figure \ref{Illustration}.
 In this regard, we define
snapshots as the spaces spanned by the offline spaces (constructed via
local spectral decomposition) at the previous step.
More precisely, at the finest coarse-grid
level, we first identify local solutions with randomized boundary conditions.
The span of these local solutions defines the snapshot space, $V_{\text{snap}}^{N-1}$.
Furthermore, we identify the local spectral problems and identify multiscale basis
functions. We call the space spanned by multiscale basis functions
the offline space and
denote the $V_{\text{off}}^{N-1}$.
These functions are assumed to represent the solution at the
$H_{N-1}$ grid. The snapshots at the $H_{N-2}$ grid can be
defined as functions, which
span all functions in the offline space, $V_{\text{off}}^{N-1}$.
However, one solves local problems in $V_{\text{off}}^{N-1}$
in an oversampled domain of a coarse region representing the coarse grid $H_{N-2}$
with the boundary conditions, which consists of randomized snapshots.
This construction constitutes a main outline of the method and a main
ingredient in homogenization methods, where the cell problems are computed
via local solutions.

The success of numerical multiscale methods depends on adaptive simulations.
The proposed method needs some adaptivity criteria to appropriately select
the number of basis function in each level. We discuss adaptivity and
present numerical results. The main idea of the adaptivity is to
start at the coarsest level $H_1$ and identify the regions, which
need additional multiscale basis functions. Then, within these regions,
we identify smaller coarse regions (at the $H_2$ level), where one needs
more basis functions. This procedure can be repeated until we reach the
$H_{N-1}$
level. In this way, we identify the regions that need additional basis
functions.

We would like to
briefly draw a parallel between the proposed methods and hierarchical
upscaling methods (such as re-iterated homogenization  and their
numerical counterparts). For this reason, we consider the flow problem
\[
\text{div} (\kappa \nabla u)=f.
\]
Our proposed approach shares a similarity with
re-iterated numerical upscaling and can be regarded
as a generalization (similar to the generalization of two-level GMsFEM and numerical upscaling,
see  \cite{chung2016_review}). In these upscaling approaches, the media properties are upscaled
at each level and then used at the next level. Numerically, this can be considered
in the following way. At the finest coarse-grid level, $H_{N-1}$, the effective properties,
$\kappa^*_{N-1}$,
are computed by solving the
local problem subject to periodic or linear boundary conditions
\[
\text{div} (\kappa \nabla \phi_i)=0 \ \text{in} \ K,
\]
$\phi_i=x_i$ on $\partial K$, for each coarse-grid block $K$ at level $N-1$.
Here, we use $K$ as a generic coarse block.
Then, the effective properties are computed in each coarse block
\[
\kappa_{K,N-1}^*={1\over |K|} \int_K \kappa \nabla \phi_i.
\]
These coarse-grid conductivities are used to solve the local
problem at the next level in the same way. The process is repeated
until we reach the desired level $\kappa^*_1$, where
the global problem is solved. As we have discussed in \cite{chung2016_review},
the limited number of local solutions (which use $x_i$ boundary conditions)
 can be thought of as multiscale basis
functions. Thus, the re-iterated homogenization or its numerical counterpart
computes the next level of multiscale basis functions using the space spanned
by the previous (finer) level multiscale basis functions. This concept is used
in our approach with the exception that computing of multiscale basis
functions is systematic and one can compute many basis functions. An important
part in the re-iterated homogenization is that one uses the same linear boundary
conditions at each level when solving local problems. Since the linear functions
can be recovered at each level exactly, one does not
need to change these boundary conditions.
In our proposed method, we use {\it all possible boundary conditions}
 (randomization is used)
spanned by multiscale basis functions of the previous (finer) coarse-grid level.
The main advantage of the proposed method is that it can handle non-separable
scales and high contrast by identifying an appropriate number of local
problems to approximate the solution at each level.

%Numerical simulations
In the paper, we consider several numerical examples. Our first set of numerical
examples is for elliptic problems in heterogeneous domains. In these examples,
we compare the results, when the multiscale basis functions are consrtucted using
a multilevel approach and when the multiscale basis functions are constructed on the
finest grid. Our numerical results show that one can achieve a similar accuracy.
We present an adaptive startegy. The main idea is to first define the coarsest region,
which needs an additional basis function. Then, within this region, we identify
subregions, which need additional multiscale basis functions. The numerical
results are presented for an adaptive method. Our second set of test examples includes
problems in perofated domains. In these examples, the domain is heterogeneous. We
present numerical results and show that one can approximate multiscale basis functions
re-iteratively.

In summary, the paper is organized as follows. In Section \ref{sec:prelim}, we present
some preliminaries and notations. In Section \ref{sec:alg}, we present the main algorithm.
Section \ref{sec:num_results} is devoted to numerical results.
In Section \ref{sec:cost}, we discuss the computational cost of the proposed method.
Finally, the paper ends with a conclusion in Section \ref{sec:conclusion}.

\section{Preliminaries}
\label{sec:prelim}

In this section, we present an overview of the main concepts
and introduce some notations.
We consider the problem
\begin{equation}
\mathcal{L} u = f \ \text{in} \ D
\end{equation}
where $\mathcal{L}$ is a differential operator, $D$ is a domain and $f$ is a given source term.
We assume that the solution $u$ satisfies typical boundary conditions,
such as the Dirichlet or the Neumann conditions on $\partial D$.
In this paper, we consider two representative cases. In the first case,
$\mathcal{L} u = \text{div} (\kappa\nabla u)$, where $\kappa$ is a heterogeneous
coefficient with multiple scales and high contrast. In the second case,
$\mathcal{L} u = \text{div} (\kappa\nabla u)$ and $D$ is a heterogeneous domain
with holes (see Figure \ref{L4-perf-mesh}).  The first case is a representative case
for problems with heterogeneous coefficients and can easily be generalized
to many other examples (\cite{chung2016_review}). The second case
represents problems in perforated domains, which occur in many applications
and can be generalized to other problems, e.g., Stokes flow in perforated domains.
Moreover, the second example shows a need for snapshot solutions.

Next, we present the coarse and the fine grids. We assume that the final
coarse-grid simulations will
be performed on a grid $\mathcal{T}_{H}^{(1)}$ with a grid size $H_1$. Our objective is
to construct multiscale basis functions on this grid.
In the two-level setting of GMsFEM (\cite{chung2016_review}), the multiscale basis functions for the grid $\mathcal{T}_H^{(1)}$
are computed by solving cell problems for each coarse region of $\mathcal{T}_H^{(1)}$.
The solutions to these cell problems are obtained numerically by using a discretization on a fine grid constructed within each coarse region.
Solving these problems can be expensive due to scale disparity
and it is therefore the purpose of this paper to investigate a re-iterated framework.
In this regard,
we assume that the grid $\mathcal{T}_H^{(1)}$ is
subdivided into a sequence of
finer grids $\mathcal{T}_H^{(l)}$ such that
$\mathcal{T}_{H}^{(l-1)}\supset\mathcal{T}_{H}^{(l)}$, $l=2,...,N$, with corresponding
grid sizes
$H_{1}>H_{2}>\cdots>H_{N}=h$. Here, $h$ is the grid size of the finest grid and $N$ is the number of grids.
Each coarse grid $\mathcal{T}_{H}^{(l)}$ consists of overlapping elements
$\omega_{j,l}$, where $j$ is the index for the $j$-th vertex in the grid $\mathcal{T}_H^{(l)}$.
For each $j$, we define $I_{j}$ as the set of all indices $m(j)$
such that
\[
\omega_{j,l-1}\supset\omega_{m(j),l}.
\]
More precisely, $\omega_{m(j),l}$
are the elements in the grid $\mathcal{T}_H^{(l)}$
that are contained in element $\omega_{j,l-1} \in \mathcal{T}_H^{(l-1)}$.
We remark that the definition of $I_j$ is dependent on the level index $l$.
Since the dependence of $I_j$ on $l$ will be clear in the context,
we will omit this dependence to simplify our notations.

Our re-iterated GMsFEM will compute the basis functions for the grid $\mathcal{T}_H^{(1)}$
using the basis functions constructed for the grid $\mathcal{T}_H^{(2)}$,
instead of a fine grid as in the two-level GMsFEM.
More generally, we will compute multiscale basis functions for the grid $\mathcal{T}_H^{(l-1)}$
by using the multiscale basis functions for the grid $\mathcal{T}_H^{(l)}$, for $l=2,3,\cdots, N$,
to approximate the cell problems.
The main idea of our re-iterated approach is that the multiscale basis functions
constructed for the level $l$ will give good approximations to the cell problems required
in finding the basis functions for the level $(l-1)$.
%The main idea of the method is to construct multiscale basis functions in each
%coarse grid level.
Next, we introduce the notations for the offline space and the snapshot space.
Assume that the multiscale basis functions are already obtained for the level $l$.
We call the space spanned by all these basis functions as the offline space,
denoted by $V_{H,\text{off}}^l$.
We notice that these basis functions are supported on the overlapping elements $\omega_{m,l}$.
For each of these overlapping elements $\omega_{m.l}$,
we denote the $j$-th basis functions by $\{\phi_{j}^{\omega_{m,l}}\}$.
The local offline space $V_{H,\text{off}}^{\omega_{m,l}}$ is the space spanned by basis functions $\{\phi_{j}^{\omega_{m,l}}\}$.
To construct the basis functions at the level $(l-1)$,
we consider an element $\omega_{k,l-1}$.
Following the idea of GMsFEM, we need a snapshot space $V_{H,\text{snap}}^{\omega_{k,l-1}} = \text{span}\{ \psi_j^{\omega_{k,l-1}} \}$
defined on $\omega_{k,l-1}$.
To construct the space $V_{H,\text{snap}}^{\omega_{k,l-1}}$,
we will solve a cell problem on $\omega_{k,l-1}$ by using the multiscale basis functions constructed
for the level $l$.
Using the snapshot space and a suitable spectral problem,
we can construct a set of offline basis $\{ \phi_j^{\omega_{k,l-1}} \}$ by selecting the dominant modes defined using the spectral problem.
The above procedure is repeated until the basis functions for the level $1$ are obtained.
The precise constructions will be given in the next section.

%We call this offline space and it will be spanned by
%$\{\phi_{j}^{\omega_{m,l}}\}$-$j^{\text{th}}$ at level
%$l$, which represents the element $\omega_{m,l}$. To construct, the offline space,
%we will use the snapshots denoted by
%$\{\psi_{j}^{\omega_{m,l}}\}$-$j^{\text{th}}$ for level
%$l$ in $\omega_{m,l}$. The offline space will be constructed via a local spectral decomposition
%in the snapshot space. We also introduce the notations for the offline and snapshot spaces
%for each $\omega_{m,l}$
%$V_{H,\text{off}}^{\omega_{m,l}}$, $V_{H,\text{snap}}^{\omega_{m,l}}$
%and also their union for the entire domain at coarse-grid level $l$,
%$V_{H,\text{off}}^{l}$,
%$V_{H,\text{snap}}^{l}$.

\subsection*{Summary of notations}

For the ease of reference, we summarize below the list of notations
defined above.

\begin{itemize}

\item $\mathcal{T}_{H}^{(l)}$ is the coarse-grid configuration at level $l$,
 $\mathcal{T}_{H}^{(l-1)}\supset\mathcal{T}_{H}^{(l)}$, $l=2,...,N$.

\item
$H_l$ is the grid size of $\mathcal{T}_H^{(l)}$ ($H_1>H_2>...>H_N$).

\item $\omega_{j,l}$ is the $j$-th coarse region at the coarse-grid level $l$.

\item $I_j$ is the set of indices $m(j)$ such that
$\omega_{j,l-1}\supset\omega_{m(j),l}$.

\item $\{\psi_{j}^{\omega_{m,l}}\}$ is the set of snapshots at  the coarse-grid level
$l$, which are supported on $\omega_{m,l}$.

\item $\{\phi_{j}^{\omega_{m,l}}\}$ is the set of offline basis at  the coarse-grid level
$l$, which are supported on $\omega_{m,l}$.

\item $V_{H,\text{snap}}^{\omega_{m,l}}$
is the local snapshot space at  the coarse-grid level $l$ for the subdomain $\omega_{m,l}$.

\item $V_{H,\text{off}}^{\omega_{m,l}}$
is the local offline space at  the coarse-grid level $l$ for the subdomain $\omega_{m,l}$.

\item $V_{H,\text{snap}}^{l}$ is the global snapshot space at  the coarse-grid level $l$.

\item $V_{H,\text{off}}^{l}$ is the global offline space at  the coarse-grid level $l$.

\end{itemize}

\section{Re-iterated GMsFEM}
\label{sec:alg}

In this section, we present the detail construction of the re-iterated GMsFEM. As we mentioned earlier that the algorithm is a
systematic extension of the two-level approach with the main idea of using offline spaces
of previous level (finer level) for construction of the snapshot space for the current level. The algorithm
is illustrated in Figure \ref{Illustration}, and we will explain next the full description of the algorithm.
%briefly describe the algorithm and
%some of its details.

The construction starts with a snapshot space for the level $(N-1)$ grid $\mathcal{T}_H^{(N-1)}$. To compute
the snapshot vectors, we will solve local problems on an oversampled region
with some appropriate randomized boundary conditions. In particular, for each oversampled region $\omega_{m,N-1}^+$, we solve
\begin{equation}
\label{eq:snap1}
\mathcal{L} \psi_j^{\omega_{m,N-1}}=0, \quad \text{in } \omega_{m,N-1}^+
\end{equation}
subject to boundary conditions $\psi_j^{\omega_{m,N-1}}=\mathcal{R}_j$,
where $\mathcal{R}_j$ is a random function, which takes an independent random
value at every boundary node. The equation (\ref{eq:snap1}) is solved on the finest grid $\mathcal{T}_H^{(N)}$ and
can therefore be thought as a discrete
equation defined on a fine grid $\mathcal{T}_H^{(N)}$.
We notice that the snapshot space $V_{H,\text{snap}}^{\omega_{m,N-1}}$
is obtained by the spans of the restrictions of the above $\psi_j^{\omega_{m,N-1}}$ in $\omega_{m,N-1}$.
We use the same notation for the solution of (\ref{eq:snap1}) and its restriction in $\omega_{m,N-1}$
to simplify the notations.
We remark that the oversampled region $\omega_{m,N-1}^+$ is usually obtained by enlarging $\omega_{m,N-1}$
by several grid blocks of the finest mesh $\mathcal{T}_H^{(N)}$.
Instead of the above construction of snapshot vectors, one can use other snapshot functions
that use all boundary conditions or all fine grid functions (see \cite{chung2016_review}).

Next, we describe the iterated procedure.
For a given level $l$, we assume that the snapshot space $V_{H,\text{snap}}^l$ is already determined.
Recall that $V_{H,\text{snap}}^{\omega_{m,l}}$ is the local snapshot space corresponding to the coarse region $\omega_{m,l}$.
We will first identify the offline space $V_{H,\text{off}}^{\omega_{m,l}}$ for this region.
To do so, we will make use of two bilinear forms $a(\cdot,\cdot)$ and $s(\cdot,\cdot)$,
and a suitable spectral problem defined in $V_{H,\text{snap}}^{\omega_{m,l}}$
to select some dominant modes.
More precisely, we consider the spectral problem of finding $\lambda\in\mathbb{R}$ and $u\in V_{H,\text{snap}}^{\omega_{m,l}}$ such that
\[
a(u,v)=\lambda s(u,v), \quad\forall v \in V_{H,\text{snap}}^{\omega_{m,l}}.
\]
We then let $\lambda_1 \leq \lambda_2 \leq \cdots \leq \lambda_K$ be the eigenvalues
and let $\widetilde{\phi}_1, \widetilde{\phi}_2, \cdots, \widetilde{\phi}_K$ be the corresponding eigenfunctions.
The dominant modes are defined as the first few eigenfunctions.
To construct the offline space $V_{H,\text{off}}^{\omega_{m,l}}$, we select the first few eigenfunctions
and define the offline basis $\phi_j^{\omega_{m,l}} = \chi_m^{(l)} \widetilde{\phi}_j$,
where $\{ \chi_m^{(l)} \}$ is the partition of unity corresponding to the grid $\mathcal{T}_H^{(l)}$.
The global offline space $V_{H,\text{off}}^{l}$ for the level $l$ is then defined as the span of all these basis functions.
We remark that the choice of the bilinear forms $a(\cdot,\cdot)$ and $s(\cdot,\cdot)$
are based on the convergence analysis.
For example, for elliptic equation considered in this paper, we will use
\begin{align*}
a(u,v) & = \int_{\omega_{m,l}} \kappa\nabla u\cdot\nabla v,\\
s(u,v) & = \int_{\omega_{m,l}} |\nabla\chi^{(l)}_m|^2\, u\,v.
\end{align*}
The reader can see \cite{chung2016_review}
for the spectral problems associated with other equations.

%Next, we define an offline space using local spectral decomposition in the snapshot space.
%In particular, we assume that two bilnear forms $a(\cdot,\cdot)$ and $s(\cdot,\cdot)$
%will be used to identify the modes and these bilinear forms are typically based on analysis.
%Then, the local spectral problem is to identify the modes  by
%Form $V_{H,\text{off}}^{\omega_{m,l-1}}=\text{span}\{\phi_{j}^{\omega_{m,l-1}}\}$
%by solving $u\in V_{H,\text{snap}}^{\omega_{m}}$
%\[
%a(u,v)=\lambda s(u,v)\;\forall v \in V_{H,\text{snap}}^{\omega_{m,l-1}}
%\]
%$\lambda_{1}\leq\lambda_{2}\leq\cdots\leq\lambda_{K}$, $\tilde{\phi}_{1},\dots\tilde{\phi}_{K}.$
%Then $\phi_{j}^{\omega_{m,l-1}}=\chi_{m}^{(l-1)}\tilde{\phi}_{j}.$
%For example, for elliptic equation, we will use
%\begin{align*}
%a(u,v) & = \int_{\omega_{m,l-1}} \kappa\nabla u\cdot\nabla v,\\
%s(u,v) & = \int_{\omega_{m,l-1}} |\nabla\chi^{(l-1)}_m|^2 uv
%\end{align*}

Once the offline space $V_{H,\text{off}}^{l}$ for the level $l$ is determined, we will construct the snapshot space
for the level $(l-1)$. The main idea is that we will use the basis functions in $V_{H,\text{off}}^{l}$
to solve the cell problems in the level $(l-1)$.
Consider a coarse region $\omega_{m,l-1}$ in the grid $\mathcal{T}_H^{(l-1)}$.
We will need to construct the snapshot space $V_{H,\text{snap}}^{\omega_{m,l-1}}$,
which is spanned by a set of basis functions $ \psi_j^{\omega_{m,l-1}}$.
Conceptually, we will solve
\begin{equation}
\label{eq:snap_l}
\mathcal{L} \psi_j^{\omega_{m,l-1}}=0, \quad \text{in } \omega_{m,l-1}.
\end{equation}
To do so, we need to specify a discretization and a boundary condition.
We recall that $I_m$ is the set of indices $k$ such that $\omega_{k,l} \subset \omega_{m,l-1}$.
For the discretization of (\ref{eq:snap_l}), we will assume that $\psi_j^{\omega_{m,l-1}}$
is a linear combination of all the offline basis $\phi_j^{\omega_{k,l}} \in V_{H,\text{off}}^{\omega_{k,l}}$
where $k\in I_m$.
For the boundary condition for (\ref{eq:snap_l}),
we will take the restriction of $\phi_j^{\omega_{k,l}}$ on the boundary of $\omega_{m,l-1}$
with the condition that the restriction of non-zero.
Thus, we have a set of linearly independent boundary conditions
and these will generate the snapshot space $V_{H,\text{snap}}^{\omega_{m,l-1}}$.
On the other hand, one can do the above computations on an oversampled region $\omega_{m,l-1}^+$,
obtained by enlarging $\omega_{m,l-1}$ by few grid blocks in $\mathcal{T}_H^{(l)}$.

The above process is repeated until the coarsest grid, that is $l=1$, where multiscale basis functions are
constructed and used for the numerical simulations.
We summarize below the main steps of the algorithm.

\subsection*{The re-iterated GMsFEM}

\noindent
\underline{Initialization}: Construct the snapshot space $V_{H,\text{snap}}^{\omega_{m,N-1}}$. Solve
\begin{equation}
\label{eq:snap1a}
\mathcal{L} \psi_j^{\omega_{m,N-1}}=0, \quad \text{in } \omega_{m,N-1}^+
\end{equation}
subject to boundary conditions $\psi_j^{\omega_{m,N-1}}=\mathcal{R}_j$,
where $\mathcal{R}_j$ is a random function. Then take the restriction of $\psi_j^{\omega_{m,N-1}}$ in $\omega_{m,N-1}$.
We perform the above computations for each coarse region $\omega_{m,N-1}$ in the grid $\mathcal{T}_H^{(N-1)}$.

\vspace{0.2cm}
\noindent
\underline{Iteration}: For each $l=N-1,\cdots, 1$, perform the following computations
\begin{itemize}
\item Construct the offline space $V_{H,\text{off}}^{\omega_{m,l}}$. Solve the spectral problem of finding $\lambda\in\mathbb{R}$ and $u\in V_{H,\text{snap}}^{\omega_{m,l}}$ such that
\[
a(u,v)=\lambda s(u,v), \quad\forall v \in V_{H,\text{snap}}^{\omega_{m,l}},
\]
and select the dominant modes.

\item Construct the snapshot space $V_{H,\text{snap}}^{\omega_{m,l-1}}$. Solve the following equation
\begin{equation}
\label{eq:snap_la}
\mathcal{L} \psi_j^{\omega_{m,l-1}}=0, \quad \text{in } \omega_{m,l-1}
\end{equation}
by using the offline spaces $V_{H,\text{off}}^{\omega_{k,l}}$, for all $k\in I_m$.
The functions $\psi_j^{\omega_{m,l-1}}$ spans the space $V_{H,\text{snap}}^{\omega_{m,l-1}}$.

\end{itemize}

\subsection{Adaptivity}

In this part, we present an adaptive procedure for the re-iterated GMsFEM.
The multilevel adaptive enrichment procedure is designed based on the two-level adaptive enrichment
algorithm developed and analyzed in \cite{chung2014adaptive}.
The main idea is to select coarse regions with larger errors,
and then enrich the multiscale approximation space
by adding basis functions in those selected coarse regions.
There are two possibilities. First, we can use new multiscale
basis functions at all pre-identified levels to update multiscale
solution. Secondly, we can use new multiscale basis functions
to update the basis functions at the coarsest level.

In the two-level setting, that is only the meshes $H_1$ and $H_N$ are used,
we will select coarse regions in the mesh $H_1$
by computing a local residual.
The local residual for the coarse region $\omega_{i,1}$ is defined as
\begin{equation}
\label{res}
R_i(v) = \int_{\omega_{i,1}} f v - \int_{\omega_{i,1}} \kappa \nabla u \cdot \nabla v, \quad\quad v \in V_{H,\text{off}}^{\omega_{i,N}}.
\end{equation}
The norm of the residual is defined as
\begin{equation}
\label{res-norm}
\| R_i \| = \sup_{v \in V_{H,\text{off}}^{\omega_{i,N}}} \frac{|R_i(v)|}{\|v\|_{\omega_{i,1}}},
\end{equation}
where $\|v\|^2_{\omega_{i,1}} = \int_{\omega_{i,1}} \kappa |\nabla v|^2$.
We remark that the residual (\ref{res}) and its norm (\ref{res-norm})
are computed in the mesh $H_N=h$, which is the finest mesh.
We define $\eta_i^2 = \|R_i\|^2 \lambda_{l_i+1}^{-1}$, and we enumerate the coarse regions in the mesh $H_1$ such that
$$
\eta_1^2 \leq \eta_2^2 \leq \cdots,
$$
where $l_i$ is the number of offline basis chosen for the region $\omega_{i,1}$.
For a given real number $0<\theta<1$, we can choose the coarse regions so that
\begin{equation}
\label{criteria}
\theta \sum_{i=1}^N \eta_i^2 < \sum_{i=1}^k \eta_i^2.
\end{equation}
For those selected regions, we will add new offline basis functions in the approximation space. Secondly, we can use these basis functions to update
the coarsest multiscale basis functions.
The above process is terminated until a certain tolerance is reached.

Now we present the our adaptive enrichment procedure for several levels.
The idea is that, to compute the residual and its norm, we will use the mesh in the next finer level.
That is, for the residual in the mesh $H_l$, we will compute it using the mesh $H_{l+1}$.
In particular, we define
local residual for the coarse region $\omega_{i,l}$ as
\begin{equation}
\label{res-multi}
R_{i,l}(v) = \int_{\omega_{i,l}} f v - \int_{\omega_{i,l}} \kappa \nabla u \cdot \nabla v, \quad\quad v \in V_{H,\text{off}}^{\omega_{i,l+1}}
\end{equation}
with its norm defined by
\begin{equation}
\label{res-norm-1}
\| R_i \| = \sup_{v \in V_{H,\text{off}}^{\omega_{i,l+1}}} \frac{|R_i(v)|}{\|v\|_{\omega_{i,l}}}
\end{equation}
where $\|v\|^2_{\omega_{i,l}} = \int_{\omega_{i,l}} \kappa |\nabla v|^2$.
This can provide a saving in computational time
compared with the two-level approach (\ref{res})-(\ref{res-norm}).

For each iteration, we will determine coarse regions in the mesh $H_1$ using the criteria defined in (\ref{criteria}).
Then we will enrich the solution space by adding basis functions from $V_{H,\text{snap}}^{\omega_{i,1}}$.
This is performed by adding the next eigenfunctions of the spectral problem.
We remark that, by merely adding basis functions from the space $V_{H,\text{snap}}^{\omega_{i,1}}$,
the solution may not be improved when a certain number of basis functions from $V_{H,\text{snap}}^{\omega_{i,1}}$ is included.
One can observe this by the relative eigenvalue decay or the residual of the solution within coarse regions.
To further enhance the accuracy of the solution, one need to include basis functions in the next level.
In particular, for the selected coarse regions in the mesh $H_1$, we will compute
the residual $R_{j,2}$ for the mesh $H_2$ for all sub coarse regions.
Then, using the criteria in (\ref{criteria}), we can select coarse regions for the mesh $H_2$
and add basis functions from the spaces $V_{H,\text{snap}}^{\omega_{j,2}}$.
We can perform this computations until the relative decay of eigenvalue is small
or the residual of the solution has little decay.
After that, we will compute the residual in the next level, mesh $H_3$,
and continue this process.

The second approach consists of using new enriched multiscale basis functions
at level $2,3,...$  for re-computing multiscale basis functions at level
$1$. In this approach, additional multiscale basis functions can be used
for re-computing and adding new multiscale basis functions at level $1$.
Alternatively, one can repeat the second approach for updating
multiscale basis functions at level $3,...$ or use the first approach
and select new multiscale basis functions at level $2,3,...$.
In the paper, we will focus on the first approach.

Below is a flow chart of the algorithm in the three coarse-grid
level setting ($N=3$).
Assume the solution $u^{(n)} \in V^{(n)}$
at the $n$-th iteration is given.
We will construct the solution $u^{(n+1)}$
by adding new basis functions in some selected coarse regions
and some selected level $l$.
To do so, we perform the following computations.

\vspace{0.2cm}

\noindent
{Step 1: \underline{Computing a solution}:
\begin{itemize}
\item Compute a solution ${u}^{(n)} \in V^{(n)}$.
\end{itemize}
}

\vspace{0.2cm}

\noindent
 {Step 2:} \underline{Adding level 1 basis}:
\begin{itemize}
\item Compute the residual norm $\|R_{i,1}\|$ in the coarse-grid level $1$ using $u^{(n)}$.
\item Select coarse regions $\omega_{i,1}$, and add level-1 basis functions.
%\item Compute a new solution $\widetilde{u}^{(n+1)}$. If it is good, then set $u^{(n+1)}=\widetilde{u}^{(n+1)}$.
%Otherwise, add level 2 basis.
\end{itemize}

\vspace{0.2cm}

\noindent
{Step 3:} \underline{Adding level 2 basis}:
\begin{itemize}
\item Compute the residual norm $\|R_{i,2}\|$ {in the regions selected in the previous step} in the coarse-grid level $2$ using $u^{(n)}$.
\item Select coarse regions $\omega_{i,2}$, and add level-2 basis functions.
%\item Compute a new solution ${u}^{(n+1)}$.
\end{itemize}

{Repeat Step 1 - 3 until the solution $u^{(n)}$ is accuracy enough. }

%\begin{enumerate}
%\item Compute the residual norm $\|R_{i,1}\|$ in the coarse-grid level $1$.

%\item Select coarse regions $\omega_{i,1}$, and add level-1 basis functions.

%\item Compute the solution again.

%\item Do one of the following:

%\begin{enumerate}
%\item if the solution in Step 3 is good, continue step 1.

%\item otherwise, compute the residual $\|R_{j,2}\|$ in the coarse-grid level $2$, and add level-2 basis functions.

%\begin{enumerate}
%\item if the solution in Step 4(b) is good, continue Step 3.
%\item otherwise, compute the residual $\|R_{k,3}\|$ in the coarse-grid level $3$, and add level-3 basis functions.
%\end{enumerate}

%\end{enumerate}

%\end{enumerate}

Notice that the above adaptive process can be easily generalized to the multi-level case.
Finally, we remark that one can also consider the goal-oriented adaptivity \cite{chung2016goal}
or the residual-driven online adaptivity \cite{chung2015residual}.

\section{Numerical results}
\label{sec:num_results}

\subsection{Numerical results for flow in heterogeneous media}

\begin{figure}[ht]
\centering
\includegraphics[scale=0.5]{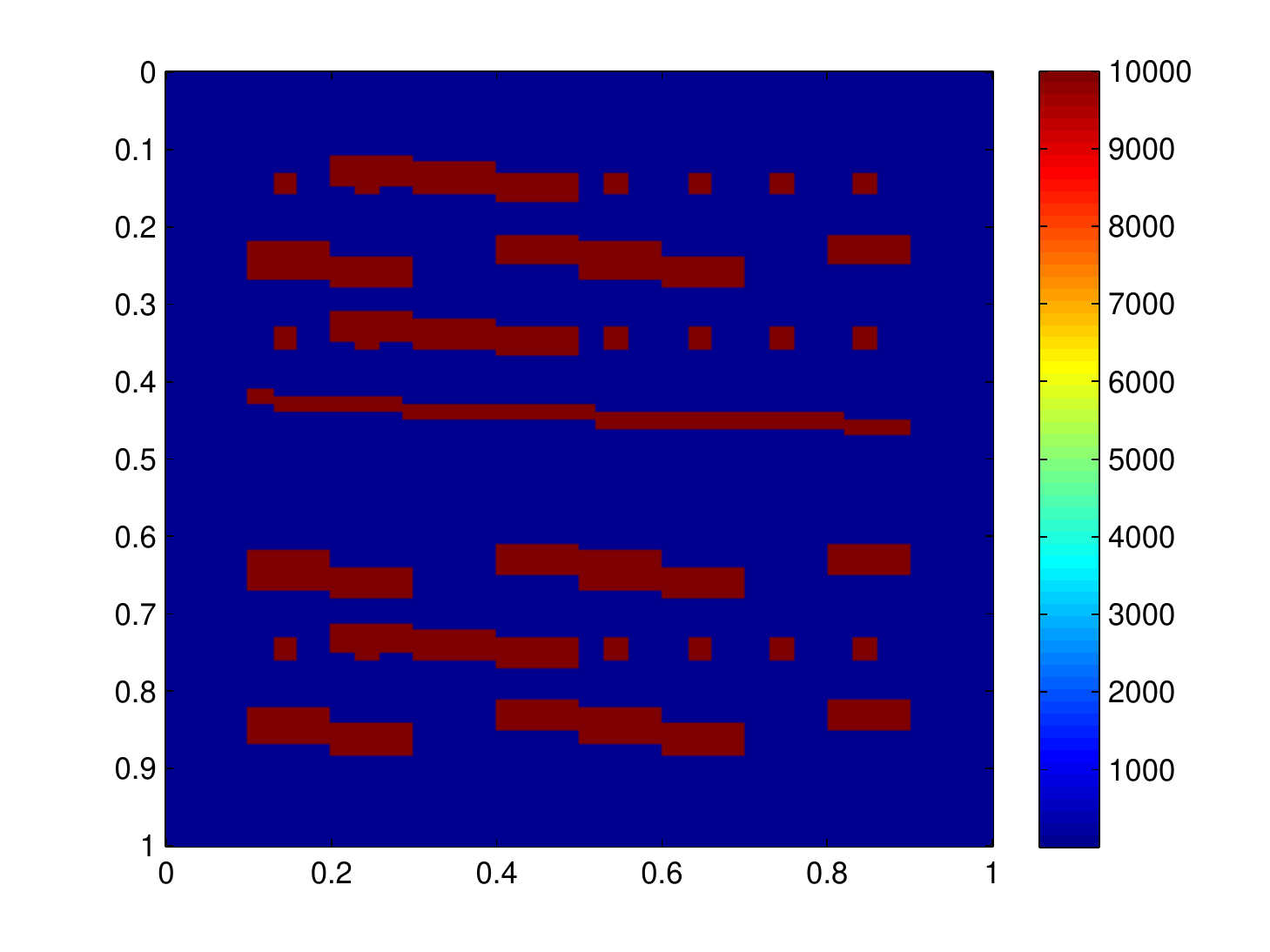}\caption{The permeability field $\kappa(x)$}
\label{fig:media}
\end{figure}

In this section, we present a number of numerical examples
to show the performance of the proposed method. The space domain $\Omega$ is taken as the unit square $[0,1]\times[0,1]$ and is divided into $10\times10$ coarse blocks consisting of uniform squares in first coarse-grid level. Each coarse block is then divided into $4\times4$ coarse blocks consisting of uniform squares in second level. Each second coarse block is divided into $10\times 10$ fine blocks consisting of uniform squares. That is, $\Omega$ is partitioned by $400\times400$ square fine-grid blocks. The high-contrast permeability field $\kappa(x)$ is shown in Figure \ref{fig:media}.

To compare the accuracy, we will use the following error quantities:
\begin{equation}\label{err_formulus}
 e_1 =\left( \frac{\|u_{H}-u_{h}\|^2_{L^2(\Omega)}}{\|u_{h}\|^2_{L^2(\Omega)}} \right)^{1/2}, \qquad
 e_2 =\left( \frac{\int_{\Omega} \kappa |\nabla(u_{H}-u_{h})|^2}{\int_{\Omega} \kappa |\nabla u_{h}|^2} \right)^{1/2}
\end{equation}
\begin{equation}\label{err_formulus2}
 e^{\text{snap}}_1 =\left( \frac{\|u_{H}-u_{\text{snap}}\|^2_{L^2(\Omega)}}{\|u_{\text{snap}}\|^2_{L^2(\Omega)}} \right)^{1/2}, \qquad
 e^{\text{snap}}_2 =\left( \frac{\int_{\Omega} \kappa |\nabla(u_{H}-u_{\text{snap}})|^2}{\int_{\Omega} \kappa |\nabla u_{\text{snap}}|^2} \right)^{1/2}
\end{equation}
where $u_h$ denotes the fine grid solution, $u_{\text{snap}}$ denotes the multiscale solution and $u_H$ denotes the multiscale solution.

%\subsubsection*{Non-adaptive result}

First, we present result for multiscale method without adaptivity. We consider the source function $f=1$. In Table \ref{tab: nonadaptive two level}, we show the convergence history for using two coarse-grid level multiscale method (the coarsest grid versus the finest grid).
We note that one can improve the results for two-level methods if mixed
or hybridization is used \cite{chung2016adaptive}. Here, our goal is to
get a comparable error to the two-level using a simplified basis construction.
In Table \ref{tab: nonadaptive three level}, we show the convergence history
for using three coarse-grid level multiscale method. For the same number of first
coarse-grid
level
basis functions, we have two choices for the second level basis
functions $10$ and $12$. We observe that the errors are similar
to those using two-level approaches. This indicates that one can use
more coarse-grid levels to inexpensively compute multiscale basis functions
and achieve similar convergence.

\begin{table}[!h]
\begin{centering}
\begin{tabular}{|c|c|c|}
\hline
$N_{1}$ & $e_{2}$ & $e_{1}$\tabularnewline
\hline
1 & 47.15\% & 28.49\%\tabularnewline
\hline
2 & 26.40\% & 7.00\%\tabularnewline
\hline
4 & 21.01\% & 4.44\%\tabularnewline
\hline
6 & 18.90\% & 3.57\%\tabularnewline
\hline
8 & 17.25\% & 2.96\%\tabularnewline
\hline
10 & 16.23\% & 2.62\%\tabularnewline
\hline
12 & 15.37\% & 2.35\%\tabularnewline
\hline
\end{tabular}
\par\end{centering}

\protect\caption{The result for two levels method. $N_1$ denotes the number of basis functions used per coarse block in first coarse-grid level.}
\label{tab: nonadaptive two level}
\end{table}

\begin{table}[!h]
\begin{centering}
\begin{tabular}{|c|c|c|c|c|c|}
\hline
$N_{1}$ & $N_{2}$ & $e_{2}$ & $e_{1}$ & $e^{\text{snap}}_{2}$ & $e^{\text{snap}}_{1}$\tabularnewline
\hline
2 & 10 & 26.06\% & 6.78\% & 20.78\% & 4.52\%\tabularnewline
\hline
4 & 10 & 21.02\% & 4.44\% & 14.00\% & 2.12\%\tabularnewline
\hline
6 & 10 & 19.10\% & 3.64\% & 11.02\% & 1.27\%\tabularnewline
\hline
8 & 10 & 17.99\% & 3.22\% & 9.02\% & 0.82\%\tabularnewline
\hline
10 & 10 & 17.61\% & 3.09\% & 8.42\% & 0.68\%\tabularnewline
\hline
12 & 10 & 17.37\% & 3.00\% & 8.01\% & 0.59\%\tabularnewline
\hline
\end{tabular} %
\begin{tabular}{|c|c|c|c|c|c|}
\hline
$N_{1}$ & $N_{2}$ & $e_{2}$ & $e_{1}$ & $e^{\text{snap}}_{2}$ & $e^{\text{snap}}_{1}$\tabularnewline
\hline
2 & 12 & 26.05\% & 6.78\% & 20.77\% & 4.52\%\tabularnewline
\hline
4 & 12 & 21.00\% & 4.44\% & 14.00\% & 2.12\%\tabularnewline
\hline
6 & 12 & 19.12\% & 3.65\% & 11.07\% & 1.28\%\tabularnewline
\hline
8 & 12 & 17.95\% & 3.20\% & 9.00\% & 0.81\%\tabularnewline
\hline
10 & 12 & 17.54\% & 3.06\% & 8.30\% & 0.65\%\tabularnewline
\hline
12 & 12 & 17.29\% & 2.97\% & 7.86\% & 0.56\%\tabularnewline
\hline
\end{tabular}
\par\end{centering}

\protect\caption{The result for three coarse-grid level method. $N_1$ denotes the number of basis functions used  per coarse block in first coarse-grid level. $N_2$ denotes the number of basis  functions used  per coarse block in second coarse-grid level.}
\label{tab: nonadaptive three level}
\end{table}

\subsection*{Adaptivity}

%$DOF_{1}:$ the number of basis used in Level 1.
%
%$DOF_{2}:$ the number of basis used in Level 2.
%
%snapshot solution is computed by using the snapshot space in level
%2.
%
%The error indicator for level 1 and level 2 are $\eta_{1}^{(i)}=\sup_{v\in V_{\text{snap},i}^{1}}\cfrac{|R(v)|^{2}}{a(v,v)},$
%$\eta_{2}^{(i)}=\sup_{v\in V_{\text{snap},i}^{2}}\cfrac{|R(v)|^{2}}{a(v,v)}$

\begin{figure}[ht]
\centering
\includegraphics[scale=0.5]{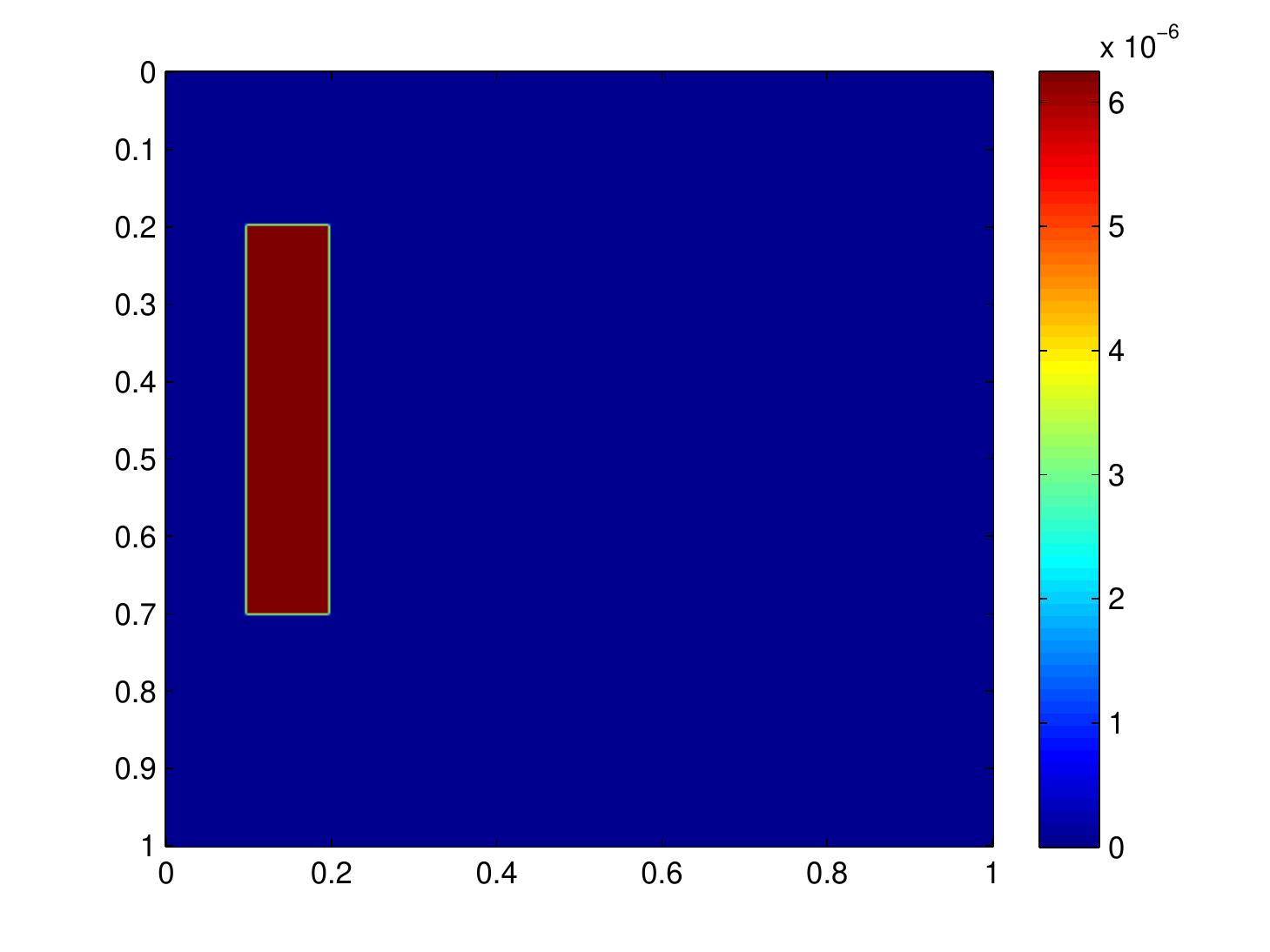}\caption{source $f$ }
\label{fig: source function}
\end{figure}

In this section, we present the numerical result for the adaptive  multiscale method. As we mentioned that one can adaptively update the number of
multiscale basis functions. This can reduce the degrees of freedom needed to represent local basis functions. In our numerical examples, we consider a different source function, $f$, which is shown in Figure \ref{fig: source function}.

\begin{figure}[ht]
\centering
\includegraphics[scale=0.5]{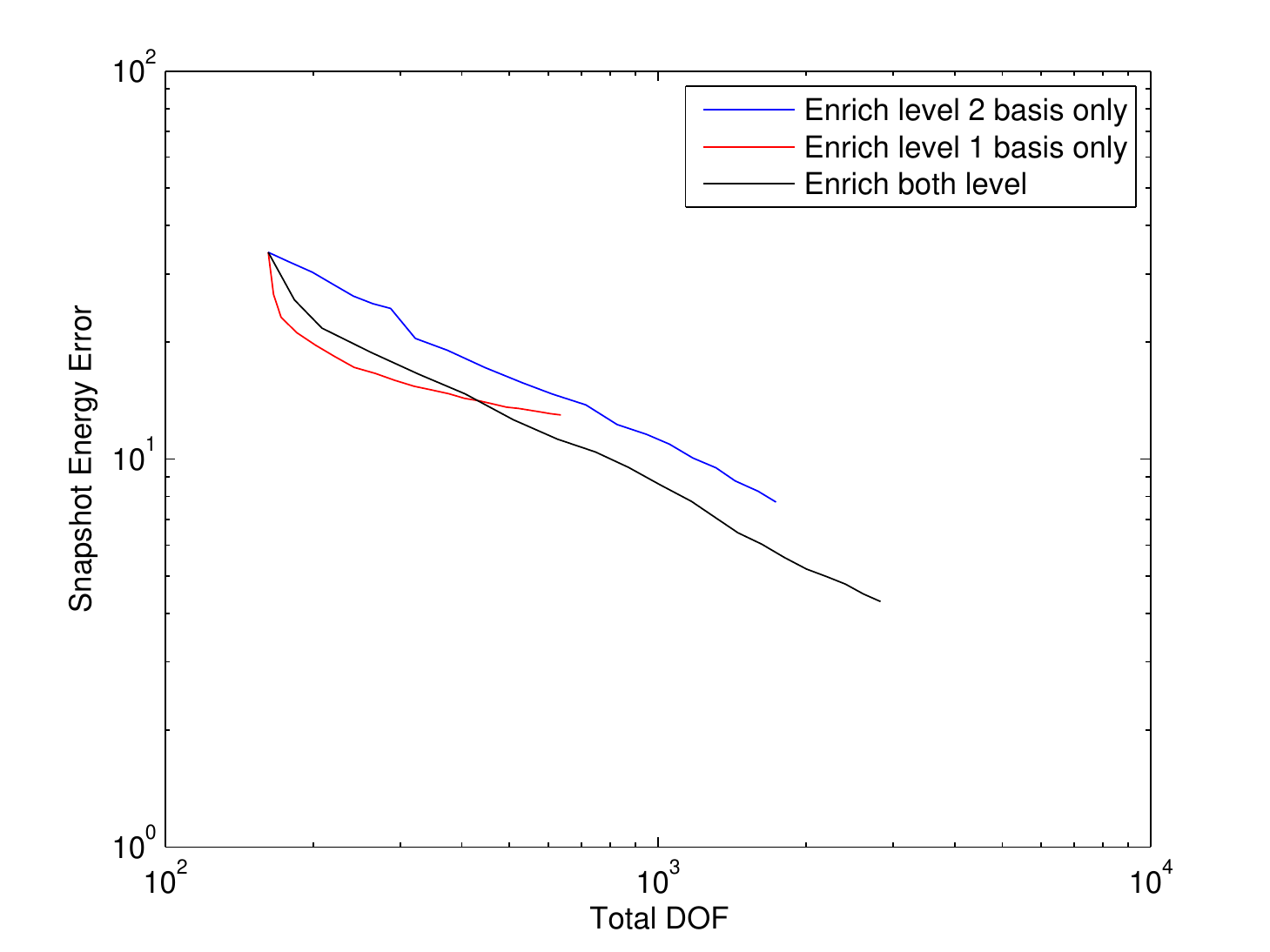}\caption{log-log scale error comparison}
\label{fig: adaptive result}
\end{figure}

In Table \ref{tab: adaptive result} and Figure \ref{fig: adaptive result}, we show the results for enriching basis in different coarse-grid levels. As we observe from this figure that one can lower the error using adaptive methods.
In particular, in the first table in Table \ref{tab: adaptive result}, we present the errors when
we adaptively enrich the basis functions in both level 1 and level 2. We observe that this gives
much better results when basis functions are only enriched in one level (results are shown in the other two tables).
More clearly, one can see in Figure \ref{fig: adaptive result} a comparison among these 3 cases,
and we observe that enriching basis in both levels gives the best results.

\begin{table}[!h]
\begin{centering}
\begin{tabular}{|c|c|c|c|}
\hline
$DOF_{1}$ & $DOF_{2}$ & $e^{\text{snap}}_{2}$ & $e^{\text{snap}}_{1}$\tabularnewline
\hline
162 & 0 & 34.12\% & 7.57\%\tabularnewline
\hline
188 & 74 & 18.82\% & 4.18\%\tabularnewline
\hline
247 & 262 & 12.63\% & 2.17\%\tabularnewline
\hline
326 & 545 & 9.51\% & 1.28\%\tabularnewline
\hline
403 & 906 & 7.06\% & 0.74\%\tabularnewline
\hline
480 & 1329 & 5.57\% & 0.46\%\tabularnewline
\hline
\end{tabular} %
\begin{tabular}{|c|c|c|c|}
\hline
$DOF_{1}$ & $DOF_{2}$ & $e^{\text{snap}}_{2}$ & $e^{\text{snap}}_{1}$\tabularnewline
\hline
162 & 0 & 34.12\% & 7.57\%\tabularnewline
\hline
185 & 0 & 21.18\% & 4.90\%\tabularnewline
\hline
242 & 0 & 17.23\% & 3.62\%\tabularnewline
\hline
320 & 0 & 15.41\% & 3.03\%\tabularnewline
\hline
405 & 0 & 14.34\% & 2.69\%\tabularnewline
\hline
493 & 0 & 13.61\% & 2.46\%\tabularnewline
\hline
\end{tabular} %
\begin{tabular}{|c|c|c|c|}
\hline
$DOF_{1}$ & $DOF_{2}$ & $e^{\text{snap}}_{2}$ & $e^{\text{snap}}_{1}$\tabularnewline
\hline
162 & 0 & 34.12\% & 7.57\%\tabularnewline
\hline
162 & 57 & 28.22\% & 6.08\%\tabularnewline
\hline
162 & 125 & 24.45\% & 5.32\%\tabularnewline
\hline
162 & 284 & 17.20\% & 3.66\%\tabularnewline
\hline
162 & 552 & 13.78\% & 2.37\%\tabularnewline
\hline
162 & 894 & 10.92\% & 1.50\%\tabularnewline
\hline
\end{tabular}
\par\end{centering}

\protect\caption{Comparison for enriching basis in different coarse-grid levels. $DOF_1$ denotes the number of basis funtions used in the first coarse-grid level. $DOF_2$ denotes the number of basis funtions used in the second coarse-grid level. }

\label{tab: adaptive result}
\end{table}

\subsection{Numerical results for perforated domain}

In this section, we consider a second example, a problem in perforated domain.
We have studied these examples for two-level approaches in our earlier
works \cite{chung2015generalized,chung2016mixed}.
We present the numerical results for four coarse-grid levels.
The computational domain $\Omega$ is taken to be $[0,1]\times [0,0.8]$ with
$230$ random  perforations. We solve elliptic equaion with source function $f=1$ with zero Dirichlet boundary conditions on perforations and global boundary.
Fine-scale mesh contains $164835$ vertices and $326048$ triangular cells.
We will use several numbers of the nested (embedded) coarse grids.
Size of the fine-scale system is $DOF_f = 164 835$.

\begin{figure}[!htb]
  \centering
  \includegraphics[width=1\textwidth]{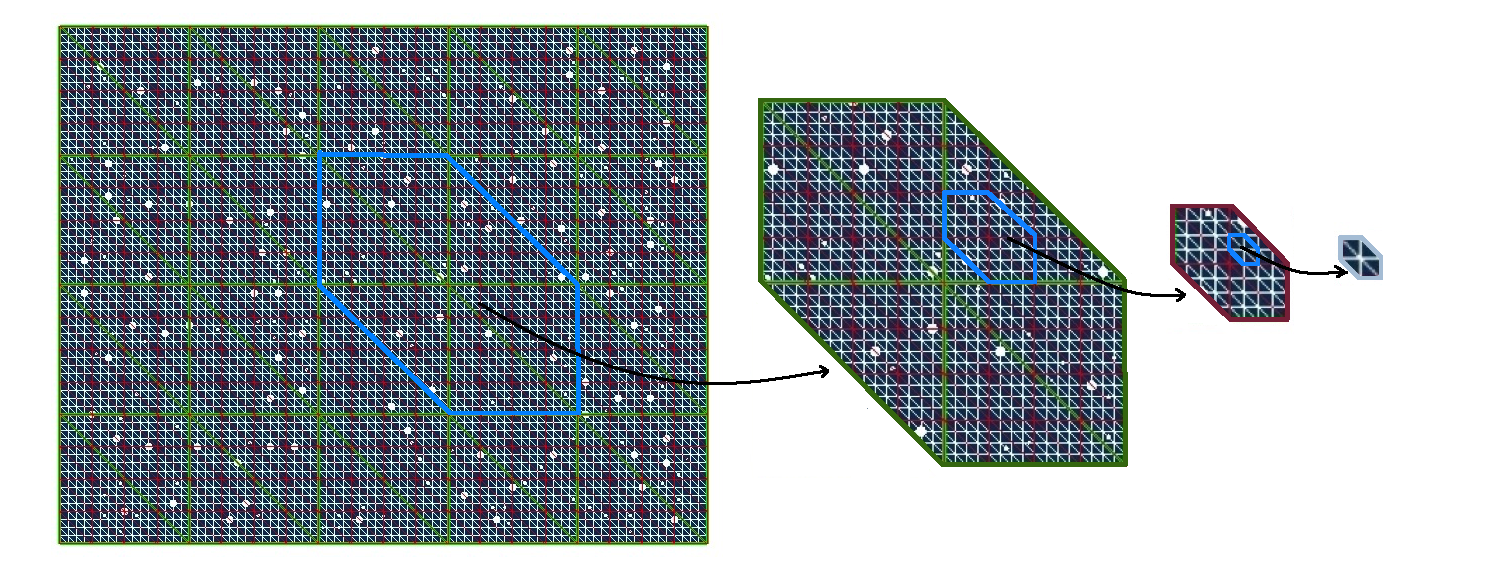}
  \caption{Computational meshes for perforated domain (230 random circle perforations) with 4 Levels. Level 1: green. Level 2: red. Level 3: white. Level 4 (fine grid): blue. Fine grid contains 164835 vertices and 326048 cells.}
 \label{L4-perf-mesh}
\end{figure}

First, we consider a case with a homogeneous background, i.e., the permeability outside perforations is $1$. We divide the fine grid into several number of coarse-grid levels and present numerical results for the following configurations:
\begin{itemize}
\item Using 2 Levels = $L1+L4$;
\item Using 3 Levels = $L1+L2+L4$;
\item Using 4 Levels = $L1+L2+L3+L4$.
\end{itemize}
Here, we use
\begin{itemize}
\item Level 1: 30 vertices and 40 cells (4x5);
\item Level 2: 357 vertices and 640 cells (4x4 for each cell from Level 1);
\item Level 3: 5265 vertices and 10240 cells (4x4 for each cell from Level 2);
\item Level 4 (fine-grid): 164835 vertices and 326048 cells.
\end{itemize}
At each level, we use a triangular grid.
For the illustration of the mesh levels, we refer to Figure \ref{L4-perf-mesh}.
Two coarse-grid level $L1+L4$ refers to the computations, where the fine grid $L4$ is
used to compute multiscale basis functions $L1$. Three coarse-grid level
$L1+L2+L4$ refers to the computation, where $L2+L4$ is used to compute
multiscale basis functions. Four coarse-grid level
$L1+L2+L3+L4$ refers to the computations, where $L2+L3+L4$
is used to compute
multiscale basis functions.

\begin{figure}[!htb]
  \centering
  \includegraphics[width=0.49\textwidth]{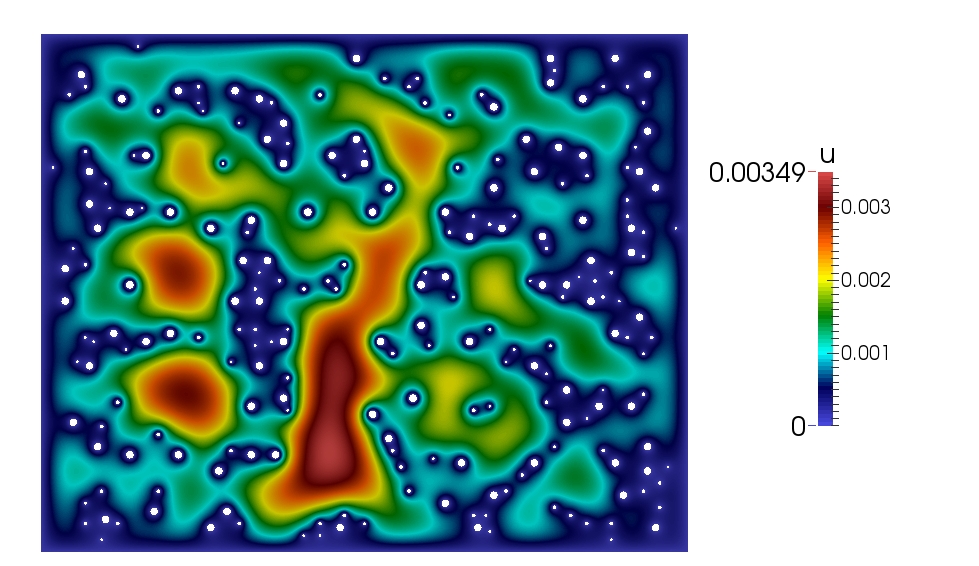}
  \includegraphics[width=0.49\textwidth]{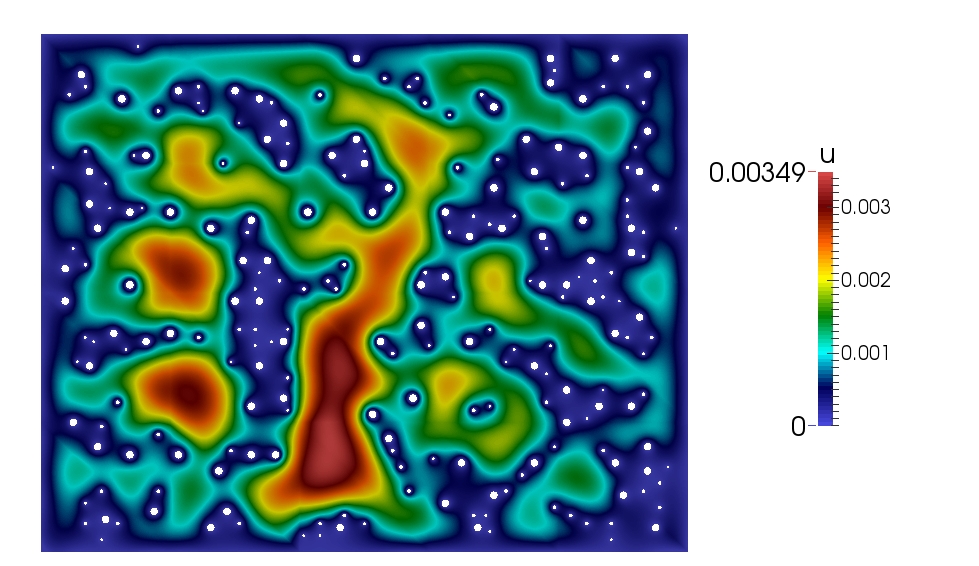}
  \caption{Fine-scale solution (left)  and coarse-scale solution (right)
 using $8$ multiscale basis functions at level $L1$ ($N_1=8$) and
$8$ basis functions for each other levels  $N_2 = N_3 = 8$.}
% 8  (bottom left) and 16  (bottom right) multiscale basis functions in Level 1 for perforated domain using 4 Levels. $N_2 = N_3 = 8$. $DOF_f = 164 835$. For $N_1 = 4$ we have $24.44 \%$ $L_2$-error,  $8.16 \%$ for $N_1 = 8$ and $3.68 \%$ for $N_1 = 16$. }
 \label{L4-perf-u}
\end{figure}

\begin{figure}[!htb]
  \centering
  \includegraphics[width=0.49\textwidth]{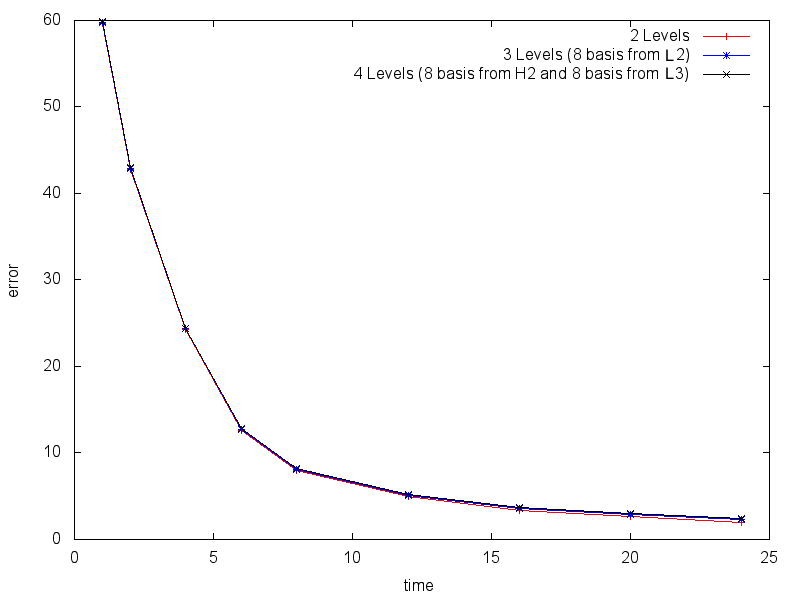}
  \includegraphics[width=0.49\textwidth]{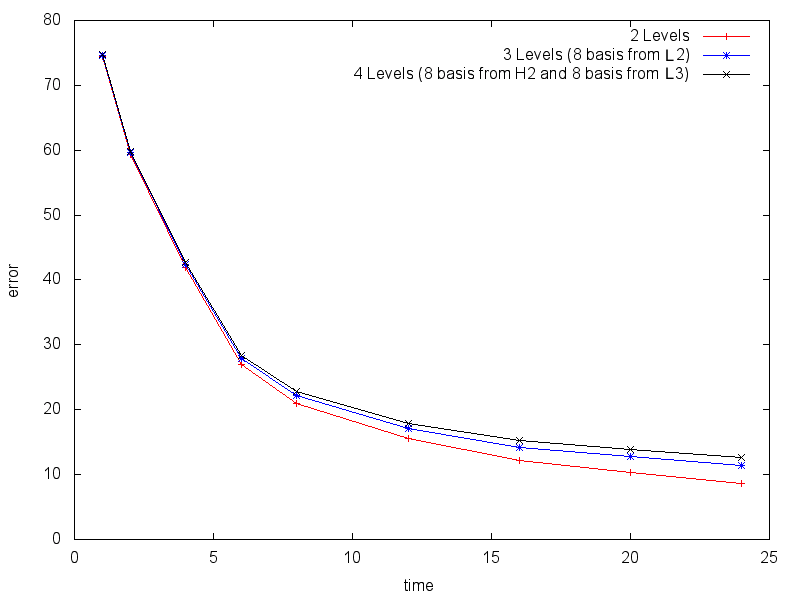}
  \caption{Illustration of Table \ref{tab-perf-dL}. Relative $L_2$( left) and energy (right) errors for perforated domain. }
 \label{L4-perf-tab1}
\end{figure}

\begin{figure}[!htb]
  \centering
  \includegraphics[width=0.49\textwidth]{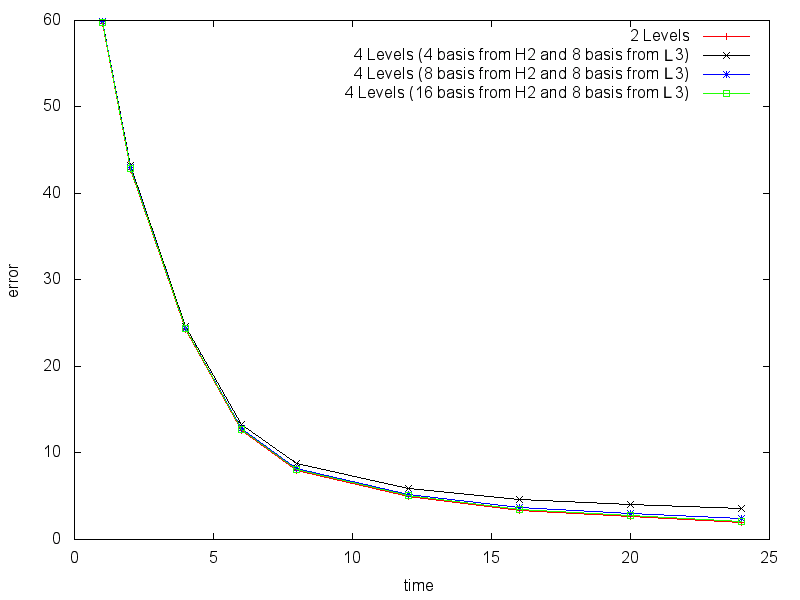}
  \includegraphics[width=0.49\textwidth]{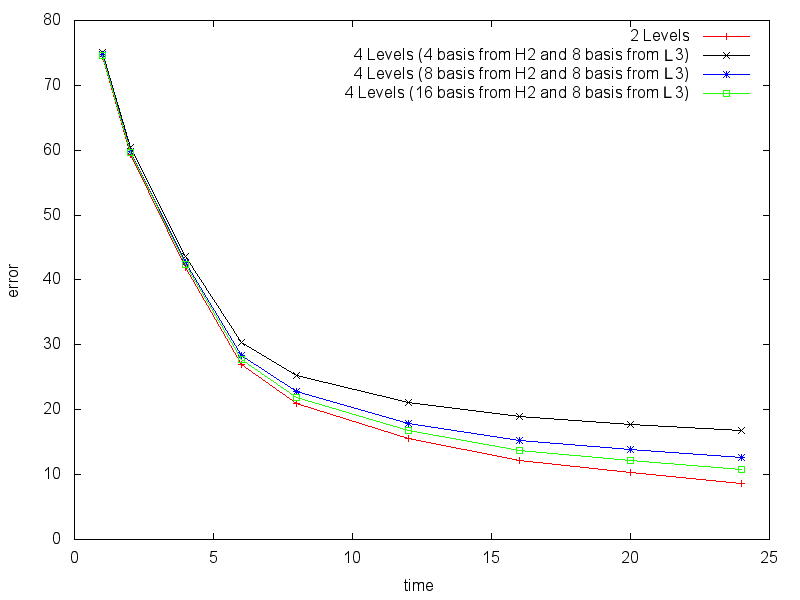}
  \caption{Illustration of Table \ref{tab-perf-L4}. Relative $L_2$ (left) and energy (right) errors for perforated domain.}
 \label{L4-perf-tab3}
\end{figure}

In Figure \ref{L4-perf-u}, we show the numerical solution of the elliptic problem using $4$ coarse-grid level multiscale method using different number of the multiscale basis at first coarsest level.
In this case, we select $8$ basis functions at the coarsest level ($N_1=8$) and also
use $8$ basis functions at each finer coarse-grid level for computing the multiscale basis functions.
The numerical results for various number of multiscale basis functions at the coarsest
level are presented in Table \ref{tab-perf-dL} (see also Figure \ref{L4-perf-tab1} for
illustration). In this table, we compare
the results obtained using two level approach, where multiscale basis functions
are computed on the fine grid to the results, when multiscale basis functions are
computed using $N_2=8$ and $N_2=N_3=8$, which use very few degrees of freedom to
approximate the basis functions. As we observe from this table that the errors
are very similar, i.e., one can use approximations for computing multiscale
basis functions. We observe similar results when fewer multiscale basis
functions are chosen at finer coarse-grid levels. In Table \ref{tab-perf-L4}, we present
numerical results by varying the number of multiscale basis functions at finer levels.
In particular, $N_2$ and $N_3$ are varied. Numerical results show that one can achieve a
similar accuracy. We depict the numerical errors in Figure \ref{L4-perf-tab3}.
%The numbers are illustrated in the captions along with the errors.
%For another meshes, we use $N_2 = N_3 = 8$ basis function for capture microstructure at each level. Note that, we should take sufficient number of the basis from previous level, otherwise we will pull errors to next level.
%When we take $N_1 = 4$ mutiscale basis functions at Level 1,
% $L_2$-error is $24.44 \%$,  $L_2$-error is $8.16 \%$ for $N_1 = 8$, and
% $L_2$-error is $3.68 \%$ for $N_1 = 16$.

% perforated
\begin{table}[!htp]
\begin{center}
\begin{tabular}{|c|c|cc|cc|cc|}
\hline
 & & \multicolumn{2}{|c|}{Using 2 Levels}
 & \multicolumn{2}{|c|}{Using 3 Levels}
 & \multicolumn{2}{|c|}{Using 4 Levels} \\
 & & \multicolumn{2}{|c|}{}
 & \multicolumn{2}{|c|}{($N_2 = 8$)}
 & \multicolumn{2}{|c|}{($N_2 = N_3 = 8$)} \\
$N_{1}$ & $DOF$ & $e_{2}$ & $e_{a}$  & $e_{2}$ & $e_{a}$  & $e_{2}$ & $e_{a}$ \\
\hline
1     & 30   	& 	59.69	& 74.52	& 	59.74	& 74.64  & 59.86	& 74.77	\\
2     & 60   	& 	42.78	& 59.37	& 	42.87	& 59.63  & 42.98	& 59.83	\\
4     & 120   & 	24.24	& 41.88	& 	24.35	& 42.38  & 24.44	& 42.69	\\
6     & 180   & 	12.64	& 26.99	& 	12.72	& 27.84  & 12.81	& 28.35	\\
8     & 240   & 	7.98		& 20.98	& 	8.07		& 22.12  & 8.16	& 22.76	\\
12   & 360   &  4.93		& 15.54	& 	5.06		& 17.13  & 5.15	& 17.95	\\
16   & 480   & 	3.40		& 12.12	& 	3.59		& 14.19  & 3.68	& 15.18	\\
\hline
\end{tabular}
\end{center}
\caption{Relative  $L_2$ and energy  errors for perforated domain using different number of coarse-grid levels (2, 3 and 4). $DOF_f = 164 835$. }
\label{tab-perf-dL}
\end{table}

\begin{table}[!htp]
\begin{center}
\begin{tabular}{|c|cc|cc|cc|}
\hline
 & \multicolumn{2}{|c|}{$N_3 = 4$}
 & \multicolumn{2}{|c|}{$N_3 = 8$}
 & \multicolumn{2}{|c|}{$N_3 = 16$} \\
$N_{1}$ & $e_{2}$ & $e_{a}$ & $e_{2}$ & $e_{a}$  & $e_{2}$ & $e_{a}$ \\
\hline
\multicolumn{7}{|c|}{$N_2 = 4$ }  \\
\hline
1     & 	60.27	& 75.40 & 60.21	& 75.20  & 60.14	& 75.07 \\
2     & 	43.68	& 61.14 & 43.42	& 60.58  & 43.36	& 60.38 \\
4     & 	25.13	& 44.83 & 24.93	& 43.92  & 24.84	& 43.53 \\
6     & 	13.84	& 32.06 & 13.26	& 30.23  & 13.18	& 29.60 \\
8     & 	9.28		& 27.30 & 8.619	& 25.08  & 8.49	& 24.27 \\
12   &  6.56		& 23.67 & 5.69	& 20.94  & 5.52	& 19.93 \\
16   & 	5.26		& 21.63 & 4.22	& 18.55  & 3.95	& 17.30 \\
\hline
\multicolumn{7}{|c|}{$N_2 = 8$ }  \\
\hline
1     & 	59.93	& 75.00 & 59.86	& 74.77  & 59.80	& 74.66 \\
2     & 	43.24	& 60.40 & 42.98	& 59.83  & 42.93	& 59.64 \\
4     & 	24.61	& 43.65 & 24.44	& 42.69  & 24.37	& 42.31 \\
6     & 	13.33	& 30.32 & 12.81	& 28.35  & 12.76	& 27.72 \\
8     & 	8.74		& 25.22 & 8.16	& 22.76  & 8.08	& 21.93 \\
12   &  5.90		& 21.12 & 5.15	& 17.95  & 5.03	& 16.82 \\
16   & 	4.61		& 18.89 & 3.68	& 15.18  & 3.49	& 13.73 \\
\hline
\end{tabular}
\end{center}
\caption{Relative $L_2$ and energy errors for perforated domain using 4 coarse-grid levels and  different number of $N_2$ and $N_3$. $DOF_f = 164 835$. }
\label{tab-perf-L4}
\end{table}

\begin{figure}[!htb]
  \centering
  \includegraphics[width=0.49\textwidth]{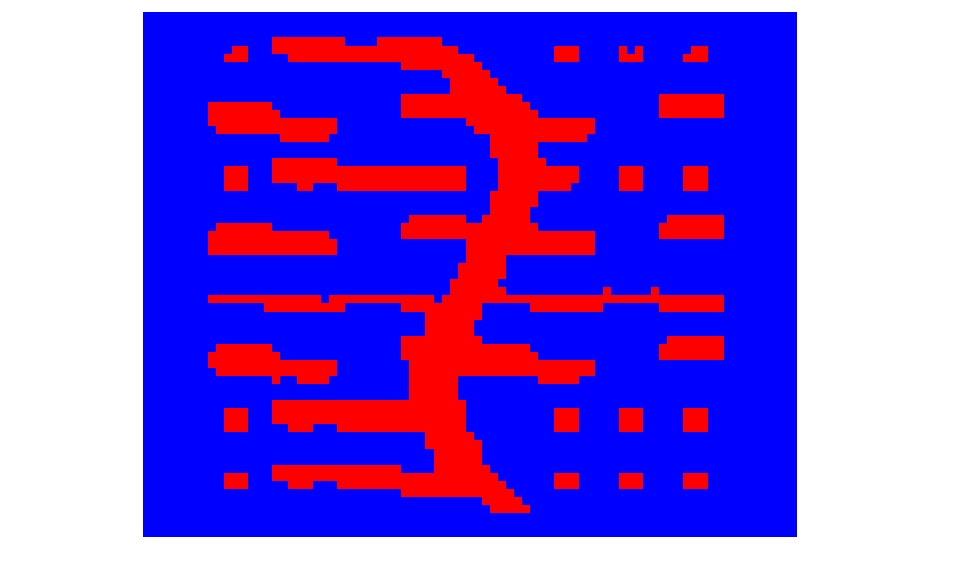}
  \includegraphics[width=0.49\textwidth]{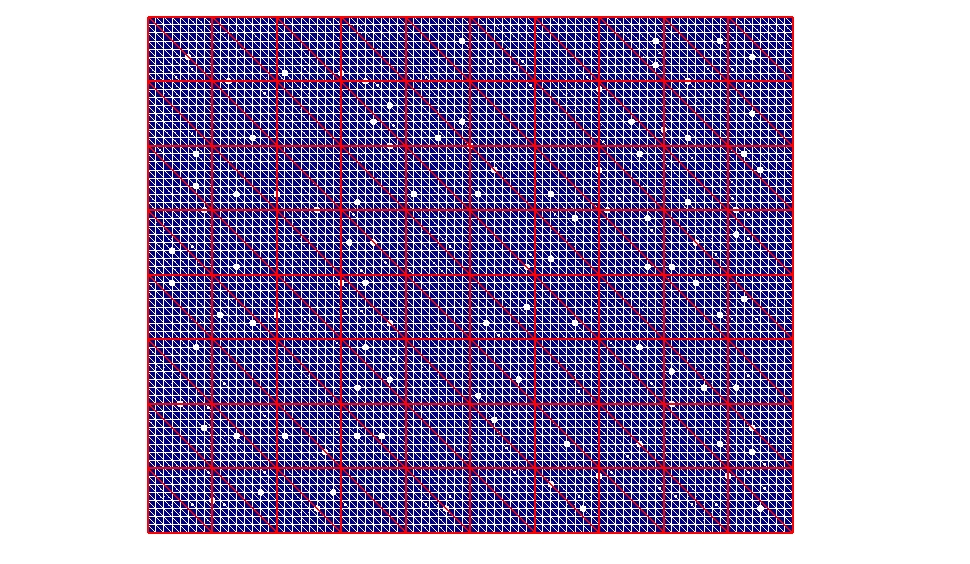}
  \caption{Heterogeneous backround (left) and computational grid (right) for perforated domain (230 random circle perforations) with three coarse-grid levels. Coarse-grid level 1: red. Coarse-grid level 2: white. Coarse-grid level 3 (fine grid): blue. Fine grid contains 164835 vertices and 326048 cells. In the left picture: $k = 1$ in red domain and $k = 1000$ in the blue domain. }
 \label{L4-heter-mesh}
\end{figure}

\begin{figure}[!htb]
  \centering
  \includegraphics[width=0.49\textwidth]{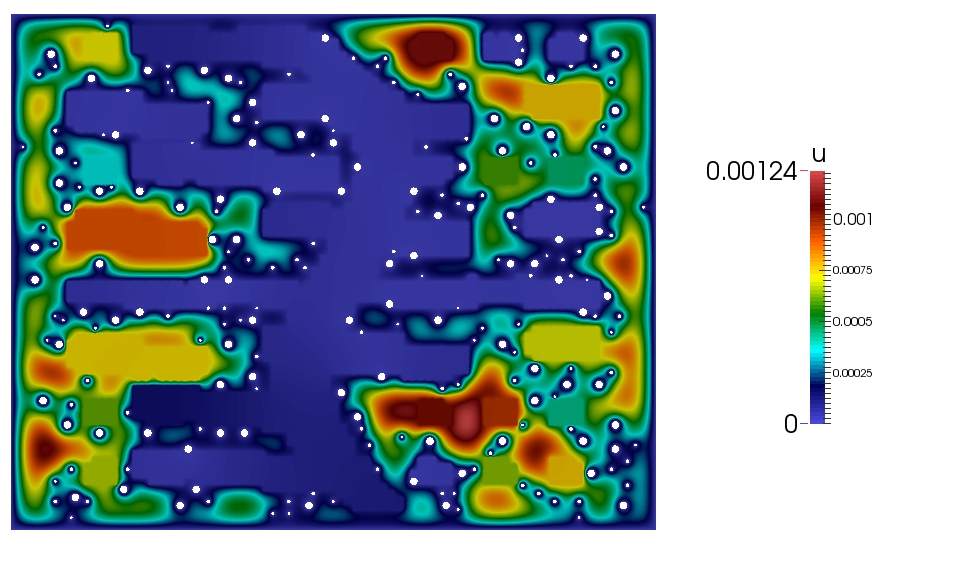}
  \includegraphics[width=0.49\textwidth]{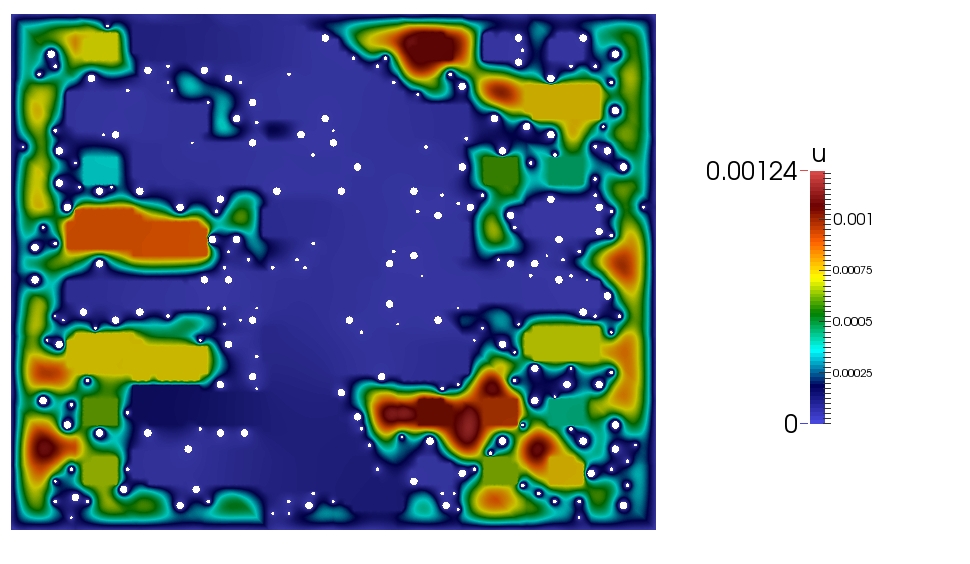}
  \caption{Fine-scale solution (top)  and coarse-scale solution using $8$ multiscale basis functions at the coarsest level and $8$ basis functions at the finer levels to compute these multiscale basis functions. $L_2$ error is $0.87$\%.}
 \label{L4-heter-u}
\end{figure}

We have observed similar results when using a heterogeneous background as
shown in Figure \ref{L4-heter-mesh} with three-level coarse grids.
\begin{itemize}
\item Level 1: 99 vertices and 160 cells (8 to 10);
\item Level 2: 1353 vertices and 2560 cells (4x4 for each cell from Level 1);
\item Level 3 (fine-grid): 164835 vertices and 326048 cells.
\end{itemize}
At each level we have triangular grid. In Figure \ref{L4-heter-u}, we depict the
fine-scale solution and the solution obtained using $8$ multiscale basis functions
at the coarsest level $N_1=8$ and $N_2=8$. In Table \ref{tab-heter-dL},
we present the errors by varying $N_1$, the number of basis functions at the coarsest level.
As we observe that approximation of multiscale basis functions gives similar
errors as two-level approaches. In this table, we also compare to the results obtained
for homogeneous background. As we observe that the errors for the homogeneous
background is better compared to those for the heterogeneous background.

% perforated
\begin{table}[!htp]
\begin{center}
\begin{tabular}{|c|c|cc|cc|}
\hline
 & & \multicolumn{2}{|c|}{Using 2 Levels}
 & \multicolumn{2}{|c|}{Using 3 Levels}\\
 & & \multicolumn{2}{|c|}{}
 & \multicolumn{2}{|c|}{($N_2 = 8$)}\\
$N_{1}$ & $DOF$ & $e_{2}$ & $e_{a}$  & $e_{2}$ & $e_{a}$ \\
\hline
\multicolumn{6}{|c|}{Heterogeneous backround}  \\
\hline
1     & 30   	& 	79.56	& 94.80	& 	79.61	& 94.81	\\
2     & 60   	& 	55.64	& 84.28	& 	55.75	& 84.36	\\
4     & 120   & 	22.36	& 57.57	& 	22.55	& 57.92	\\
6     & 180   & 	1.781	& 42.42	& 	2.11		& 43.00	\\
8     & 240   & 	1.31		& 33.98	& 	1.69		& 34.77	\\
12   & 360   &  0.67		& 23.96	& 	1.05		& 25.15	\\
16   & 480   & 	0.42		& 17.98	& 	0.87		& 19.71	\\
\hline
\multicolumn{6}{|c|}{Homogeneous backround}  \\
\hline
1     & 30   	& 	31.59	& 54.86	& 	31.86	& 55.23	\\
2     & 60   	& 	17.02	& 37.03	& 	17.26	& 37.66	\\
4     & 120   & 	7.68		& 22.70	& 	7.86		& 23.78	\\
6     & 180   & 	4.18		& 15.76	& 	4.35		& 17.31	\\
8     & 240   &  2.83		& 12.29	& 	3.01		& 14.29	\\
12   & 360   &  1.42		& 7.95		& 	1.66		& 10.84	\\
16   & 480   & 	0.89		& 5.96		& 	1.18		& 9.48		\\
\hline
\end{tabular}
\end{center}
\caption{Relative  $L_2$ and energy  errors for perforated domain with heterogeneous and homogeneous backrounds for different number of coarse-grid levels (2 and 3) $DOF_f = 164 835$. }
\label{tab-heter-dL}
\end{table}

\section{Computational cost}
\label{sec:cost}
One of our aims is to compute multiscale basis functions at the coarsest
level for problems when disparity of scales does not allow such
computations. In this sense, our approach has a similar concept
as re-iterated homogenization with the aim of computing
a reduced-order model.
For this reason, our coarsening factors are usually very large.
Besides allowing basis computations in extreme
disparate scale cases,
re-iterated approach can provide computational savings.
Below, we discuss the computational savings.

%The offline cost computations.
%The idea is to be able to compute.
%We use a very large coarsening factor.

%
%
%$N_L$ is the dimension of offline in $\omega^{+}$
%
%$N_c$ is the number of basis we will choose. Assume that we can use randomization
%
%If $N_c << N_L$, then algorithm is efficient
%
%Homogenization, we have $N_c=2$.
%
%A hierarchical representation of basis functions.

Let $C_{i}$ be the coarsening number for level $i$ with
the total number of
fine grids $N=\Pi_{i}C_{i}$.
Let $L:\mathbb{N}\rightarrow\mathbb{N}$ be the operation counting
operator for the solver of the local snapshot problem where $L(N)$
is the number of operation for solving a local problem with size $N$.
Furthermore,
let $P:\mathbb{N}\times\mathbb{N}\rightarrow\mathbb{N}$ be the operation
counting operator for the solver of the local eigenvalue problem where
$P(a,b)$ is the number of operation for finding the first $b$ eigenvector
with the local eigenvalue problem size $a$.

The operation number for constructing two level multiscale basis,
$O_{2}$, (ignoring the time for forming matrix and gridding) is
\[
O_{2}=O(C_{1}r_{1}L(M_{1}\Pi_{i>1}C_{i})+P(r_{1},\lambda_{1})),
\]
 where $M_{\ensuremath{1}}$ depends on the coarse grid structure and
the oversampling size, $r_{1}$ depends on the number of random basis
solved for each node and $\lambda_{1}$ depends on the number of basis
chosen for each node.
The operation count for constructing re-iterated multiscale basis functions,
$O_{2}$, (ignoring the time for forming matrix and griding) is
\[
O_{M}=O\left(\sum_{j=1}^{M-1}\left((\Pi_{i\leq j}C_{i})r_{i}L(\lambda_{i+1}M_{i}C_{i+1})+P(r_{i},\lambda_{i})\right)\right).
\]

We assume $M_{i}C_{i+1}>r_{i}>\lambda_{i}$,
$P(r_{i},\lambda_{i})=O(\lambda_{i}r_{i}^{\beta})$,
and $L(N)=O(N^{\alpha})$ and $\alpha\geq\beta$.
{We remark that the operator, $L$, is dependent on the contrast of the medium parameter. For a high contrast medium,  we expect to have a better performance of the local solver for re-iterated case.}
Then
\begin{align*}
O_{2} & =O(C_{1}r_{1}L(M_{1}\Pi_{i>1}C_{i}))\\
O_{M} & =O\left(\sum_{j=1}^{I-1}(\Pi_{i\leq j}C_{i})r_{i}L(\lambda_{i+1}M_{i}C_{i+1})\right).
\end{align*}
We further assume $r_{i}=r,\; M_{i}=M,\; C_{i}=C\;\forall i\geq1$ and $\lambda_{i}=\lambda\;\forall i\leq M-1$,
\begin{align*}
O_{M} & =O\left(rM^{\alpha}\left(\sum_{j=1}^{M-2}\lambda^{\alpha}(\Pi_{i\leq j}C)C^{\alpha})+C^{M-1}C^{\alpha}\right)\right)\\
 & =O\left(rM^{\alpha}C^{\alpha}C\left(\lambda^{\alpha}\cfrac{C^{M-2}-1}{C-1}+C^{M-2}\right)\right)
\end{align*}
and
\[
O_{2}=O\left(rM^{\alpha}CC^{(M-1)\alpha}\right).
\]
If $C>\lambda^{\alpha},$ we have
\begin{align*}
O_{2} & =O\left( \cfrac{C^{(M-2)\alpha}}{C^{M-2}} O_{M}\right)\\
 & =O\left(C^{(M-2)(\alpha-1)} O_{M}\right).
\end{align*}
Thus, if the coarsening factor is large and fewer basis functions
are selected (cf., re-iterated homogenization, when the number
of basis functions is very low), the re-iterated construction
will provide a computational advantage. One of additional main advantage
will be due to adaptivity, where only a few basis functions
are needed in many regions.

\section{Conclusions}
\label{sec:conclusion}

In applications, the coarse computational grid may have very high resolution
so that one can not afford local simulations to extract important modes within GMsFEM.
In these cases, approximate calculations are needed.
This paper presents an extension of two-level GMsFEM by re-iteration (similar
to homogenization)
to approximate multiscale basis functions on the coarsest level inexpensively.
Following the general concept of the GMsFEM, we define snapshots at each level,
which consist of all possible basis functions at the previous (finer) levels.
This is a natural extension and shares some similarities with re-iterated homogenization
methods, where the homogenization is carried out at separate levels and used in the
next coarse level. We present an adaptive procedure and demonstrate numerical results for
two applications.

{\bf Acknowledgements.}
EC's research was in part supported by Hong Kong RGC General Research Fund Project No. 400813. MV's  work is partially supported by Russian Science Foundation Grant RS 15-11-10024 and  RFBR 15-31-20856.

\bibliographystyle{plain}
\bibliography{references,references1}
%\bibliography{references,references1,Fractures,rlg,seismics,SeismicWaves,SeismicAnisotropy,SeismicProcessing,Inversino,refs,chester,YE_2015_references,mlmc}

\end{document}